\documentclass[a4paper,12pt]{amsart}
\usepackage{amssymb}
\usepackage{ifthen}
\usepackage[dvips]{graphicx}
\nonstopmode \numberwithin{equation}{section}
\usepackage{amssymb}
\usepackage{ifthen}
\usepackage{graphicx}
\usepackage{amsmath}
\usepackage[T1]{fontenc} 
\usepackage[utf8]{inputenc}
\usepackage[usenames,dvipsnames]{color}
\usepackage{color}
\usepackage[english]{babel}
\usepackage{fancyhdr}
\usepackage{fancybox}
\usepackage{tikz}

\nonstopmode \numberwithin{equation}{section}
\theoremstyle{plain}
\newtheorem{thm}[equation]{Theorem}
\newtheorem{cor}[equation]{Corollary}
\newtheorem{lem}[equation]{Lemma}
\newtheorem{prop}{Proposition}

\newtheorem{conj}{Conjecture}

\theoremstyle{definition}
\newtheorem{defn}{Definition}[section]

\newtheorem{prob}{Problem}
\newtheorem{rem}{Remark}[section]

\newcounter{minutes}
\divide\time by 60
\newcounter{hours}
\multiply\time by 60
\addtocounter{minutes}{-\time}

\newcounter {own}
\def\theown {\thesection       .\arabic{own}}

\newenvironment{pf}[1][]{%
	\vskip 3mm
	\noindent
	\ifthenelse{\equal{#1}{}}%
	{{\slshape Proof. }}%
	{{\slshape #1.} }%
}%
{\qed\bigskip}

\newcounter{alphabet}

\newcommand{\real}{{\operatorname{Re}\,}}

\def\be{\begin{equation}}
\def\ee{\end{equation}}

\newcommand{\bee}{\begin{enumerate}}
	\newcommand{\eee}{\end{enumerate}}

\newcommand{\blem}{\begin{lem}}
	\newcommand{\elem}{\end{lem}}
\newcommand{\bthm}{\begin{thm}}
	\newcommand{\ethm}{\end{thm}}
\newcommand{\bcor}{\begin{cor}}
	\newcommand{\ecor}{\end{cor}}
\newcommand{\beg}{\begin{examp}}
	\newcommand{\eeg}{\end{examp}}
\newcommand{\begs}{\begin{examples}}
	\newcommand{\eegs}{\end{examples}}

\newcommand{\bdefn}{\begin{defn}}
	\newcommand{\edefn}{\end{defn}}

\newcommand{\bprob}{\begin{prob}}
	\newcommand{\eprob}{\end{prob}}
\newcommand{\bei}{\begin{itemize}}
	\newcommand{\eei}{\end{itemize}}

\newcommand{\bcon}{\begin{conj}}
	\newcommand{\econ}{\end{conj}}
\newcommand{\bcons}{\begin{conjs}}
	\newcommand{\econs}{\end{conjs}}
\newcommand{\bprop}{\begin{prop}}
	\newcommand{\eprop}{\end{prop}}
\newcommand{\br}{\begin{rem}}
	\newcommand{\er}{\end{rem}}
\newcommand{\brs}{\begin{rems}}
	\newcommand{\ers}{\end{rems}}
\newcommand{\bo}{\begin{obser}}
	\newcommand{\eo}{\end{obser}}
\newcommand{\bos}{\begin{obsers}}
	\newcommand{\eos}{\end{obsers}}
\newcommand{\bpf}{\begin{pf}}
	\newcommand{\epf}{\end{pf}}
\newcommand{\ba}{\begin{array}}
	\newcommand{\ea}{\end{array}}
\newcommand{\beq}{\begin{eqnarray}}
\newcommand{\beqq}{\begin{eqnarray*}}
\newcommand{\eeq}{\end{eqnarray}}
\newcommand{\eeqq}{\end{eqnarray*}}

\begin{document}

\title{The Bohr inequality for certain harmonic mappings}

\author{Vasudevarao Allu}
\address{Vasudevarao Allu,
	School of Basic Science,
	Indian Institute of Technology Bhubaneswar,
	Bhubaneswar-752050, Odisha, India.}
\email{avrao@iitbbs.ac.in}

\author{Himadri Halder}
\address{Himadri Halder,
	School of Basic Science,
	Indian Institute of Technology Bhubaneswar,
	Bhubaneswar-752050, Odisha, India.}
\email{hh11@iitbbs.ac.in}

\subjclass[{AMS} Subject Classification:]{Primary 30C45, 30C50, 30C80}

\def\thefootnote{}
\footnotetext{ {\tiny File:~\jobname.tex,
		printed: \number\year-\number\month-\number\day,
		\thehours.\ifnum\theminutes<10{0}\fi\theminutes }
} \makeatletter\def\thefootnote{\@arabic\c@footnote}\makeatother

\begin{abstract}
Let  $\phi$ be analytic and univalent ({\it i.e.,} one-to-one) in $\mathbb{D}:=\{z\in\mathbb{C}: |z|<1\}$ 
such that $\phi(\mathbb{D})$ has positive real part, is symmetric with respect to the real axis, 
starlike with respect to $\phi(0)=1,$ and $\phi ' (0)>0$.
A function $f \in \mathcal{C}(\phi)$ if $1+ zf''(z)/f'(z) \prec \phi (z),$ and $f\in \mathcal{C}_{c}(\phi)$ if $2(zf'(z))'/(f(z)+\overline{f(\bar{z})})' \prec \phi (z)$ for $ z\in \mathbb{D}$.	In this article, we consider the classes $\mathcal{HC}(\phi)$ and $\mathcal{HC}_{c}(\phi)$ consisting of harmonic mappings $f=h+\overline{g}$ of the form 
$$
h(z)=z+ \sum \limits_{n=2}^{\infty} a_{n}z^{n} \quad \mbox{and} \quad  g(z)=\sum \limits_{n=2}^{\infty} b_{n}z^{n}
$$
in the unit disk $\mathbb{D}$, where $h$ belongs to $\mathcal{C}(\phi)$ and $\mathcal{C}_{c}(\phi)$ respectively, with the dilation $g'(z)=\alpha z h'(z)$ and $|\alpha|<1$. Using the Bohr phenomenon for subordination classes \cite[Lemma 1]{bhowmik-2018}, we find the radius $R_{f}<1$ such that Bohr inequality 
$$
|z|+\sum_{n=2}^{\infty} (|a_{n}|+|b_{n}|)|z|^{n} \leq d(f(0),\partial f(\mathbb{D}))
$$
holds for $|z|=r\leq R_{f}$ for the classes $\mathcal{HC}(\phi)$ and $\mathcal{HC}_{c}(\phi)$ . As a consequence of these results, we obtain several interesting corollaries on Bohr inequality for the aforesaid classes.
\end{abstract}

\maketitle
\pagestyle{myheadings}
\markboth{Vasudevarao Allu and  Himadri Halder}{The Bohr inequality for certain harmonic mappings}

\section{Introduction and Preliminaries}

\vspace{4mm}
In recent years studying the Bohr inequality has become an interesting topic in geometric function theory, which (in the final form has been independently proved by Weiner, Riesz and Schur) states that if $f$ is an analytic function in the unit disk $\mathbb{D}:=\{z\in \mathbb{C}: |z|<1\}$ with the following Taylor series expansion
\begin{equation} \label{him-p4-e-1.1}
f(z)=\sum_{n=0}^{\infty} a_{n}z^{n}
\end{equation}
such that $|f(z)|<1$ in $\mathbb{D}$, then the majorant series $M_{f}(r)$ associated with $f$ satisfies the following inequality
\begin{equation} \label{him-p4-e-1.2}
M_{f}(r):= \sum_{n=0}^{\infty} |a_{n}|r^{n} \leq 1 \quad \mbox{for} \quad |z|=r \leq 1/3,
\end{equation}
and the constant $1/3$, known as Bohr radius, cannot be improved further.  Analytic functions $f$ of the form \eqref{him-p4-e-1.1} with $|f(z)|<1$ satisfying the inequality \eqref{him-p4-e-1.2} for $|z|=r\leq1/3$, are sometimes  said to satisfy the classical Bohr phenomenon. It is worth noting that the inequality \eqref{him-p4-e-1.2}, called Bohr inequality, can be written in the following form 
\begin{equation} \label{him-p4-e-1.3}
\sum_{n=1}^{\infty} |a_{n}z^{n}|\leq 1-|a_{0}|=d(f(0),\partial f(\mathbb{D}))
\end{equation}
for $|z|=r \leq 1/3$ and the constant $1/3$ is independent on the coefficients of the Taylor series \eqref{him-p4-e-1.1}, where $d$ is the Euclidean distance. In a better way we can demonstrate this fact by saying that Bohr phenomenon occurs in the class of analytic self maps of the unit disk $\mathbb{D}$. In view of \eqref{him-p4-e-1.3}, the notion of the Bohr phenomenon can be generalized to the class $\mathcal{G}$ consisting of analytic functions $f$ in $\mathbb{D}$ which take values in a given domain $D \subseteq \mathbb{C}$ such that $f(\mathbb{D}) \subseteq D$ and the class $\mathcal{G}$ is said to satisfy the Bohr phenomenon if there exists largest radius $r_{D} \in (0,1)$ such that the inequality \eqref{him-p4-e-1.3} holds for $|z|=r\leq r_{D}$ and for all functions $f\in \mathcal{G}$. The largest radius $r_{D}$ is called the Bohr radius for the class $\mathcal{G}$. When $D$ is a proper simply connected domain in $\mathcal{C}$,  Abu-Muhanna \cite{Abu-2010} has proved that the sharp radius $r_{D}$ is $3-2\sqrt{2}$ for the class $\mathcal{G}$. For any convex domain $D$, Aizenberg \cite{aizn-2007} has shown that $r_{D}$ coincides with the classical Bohr radius $1/3$. The Bohr inequality has been generalized to several complex variables by finding the multidimensional Bohr radius (see \cite{aizn-2000,aizn-2007,boas-1997}). For more intriguing aspects of Bohr inequality, we refer the reader to the articles \cite{abu-2014,Ali-2017,alkhaleefah-2019,Himadri-Vasu-P1,Himadri-Vasu-P2,Himadri-Vasu-P3,kaptanoglu-2006,kayumov-2018-b,Lasere-2013}.
\vspace{3mm}

In general, one can obtain the Bohr radius for certain classes of analytic and harmonic functions in $\mathbb{D}$, when the sharp coefficient bounds for these classes are known. Indeed, the proof of the inequality \eqref{him-p4-e-1.2} follows from the sharp coefficient estimates $|a_n|\leq 1-|a_0|^2$ for $n\geq 1$, which may be obtained as an application of Pick's invariant form of Schwarz's lemma for functions $f \in \mathcal{B}$, where $\mathcal{B}$ denotes the class of analytic functions in $\mathbb{D}$ with modulus less than $1$. For $f\in \mathcal{B}$, Pick's conformally invariant formula gives:
$$
|f'(z)|\leq \frac{1-|f(z)|^2}{1-|z|^2} \quad \mbox{for} \quad z \in \mathbb{D}.
$$
Furthermore, $|f'(0)|=|a_1|\leq 1-|a_0|^2$, from which for $n \geq 1$, the sharp inequality $|a_n|\leq 1-|a_0|^2$ follows. In a better way, we can demonstrate the fact that Bohr radius is estimated for a given class of functions provided the sharp coefficient bounds of the functions in that class are known. It is worth mentioning that the sharp coefficient bounds for most of the classes of harmonic functions whose analytic parts are Ma-Minda convex functions and convex functions with respect to symmetric points, which have been considered in this paper, are not yet known. 
This shows that solving Bohr radius for these classes is an open challenging problem. In this paper, we take the opportunity to solve that problem to find the best possible lower bound of the radius so that Bohr phenomenon holds for these harmonic classes.
\vspace{4mm}

Let $\mathcal{A}$ denote the class of analytic functions in $\mathbb{D}$ with $f(0)=0$ and $f'(0)=1$. Each function $f \in \mathcal{A}$ has the following representation
\begin{equation} \label{him-p4-e-1.4}
f(z)=z+ \sum_{n=2}^{\infty} a_{n}z^{n}.
\end{equation}  
Let $\mathcal{S}$ be the subclass of $\mathcal{A}$ consisting of univalent functions.
An analytic function $f\in \mathcal{S}$ is said to be starlike (respectively convex) if
$\real\left(zf'(z)/f(z)\right)>0$ (respectively $\real \left(1+ zf''(z)/f'(z)\right)>0$) for $z\in \mathbb{D}$. Let $\mathcal{S}^{*}$ (respectively $\mathcal{C}$) be the subclass of $\mathcal{S}$ consisting of starlike (respectively convex) functions in $\mathbb{D}$. For more properties of starlike and convex functions, we suggest the reader to the book \cite{vasu-book}. Using the notion of subordination, Ma and Minda \cite{ma minda-1992-a} have introduced more general subclasses of $\mathcal{S}^{*}$ and $\mathcal{C}$, denoted by $\mathcal{S}^{*}(\phi)$ and $\mathcal{C}(\phi)$, consisting of functions in $\mathcal{S}$ for which $zf'(z)/f(z) \prec \phi (z)$ and $1+ zf''(z)/f'(z) \prec \phi (z)$ respectively. Here the function $\phi :\mathbb{D} \rightarrow \mathbb{C}$, called Ma-Minda function, is analytic and univalent in $\mathbb{D}$ such that $\phi(\mathbb{D})$ has positive real part, symmetric with respect to the real axis, starlike with respect to $\phi(0)=1$ and $\phi ' (0)>0$. A Ma-Minda function has the series representation of the form $\phi(z)=1+ \sum_{n=1}^{\infty} B_{n}z^{n} \quad (B_{1}>0)$. Similarly, it is natural to consider the function $\psi$, called non-Ma-Minda function, with the condition $\psi'(0)<0$ and the other conditions on $\psi$ are the same as that of $\phi$. Note that $\psi$ can be obtained from $\phi$ by a rotation, namely, $z$ by $-z$. Similar to the definition of $\mathcal{S}^{*}(\phi)$ and $\mathcal{C}(\phi)$, we consider the classes $\mathcal{S}^{*}(\psi)$ and $\mathcal{C}(\psi)$, where $\psi$ is non-Ma-Minda function. The extremal functions $K$ and $H$ respectively for the classes $\mathcal{C}(\phi)$ and $\mathcal{S}^{*}(\phi)$ are given by
\begin{equation} \label{him-p3-e-1.6}
1+ \dfrac{zK''(z)}{K'(z)} = \phi (z) \quad \mbox{and} \quad \dfrac{zH'(z)}{H(z)} =\phi (z)
\end{equation}
with the normalizations $K(0)=K'(0)-1=0$ and $H(0)=H'(0)-1=0$. The  functions $K$ and $H$ belong to the classes $\mathcal{C}(\phi)$ and $\mathcal{S}^{*}(\phi)$ and they play the role of the Koebe function in the respective classes. We have the following subordination theorems and growth estimates for the class $\mathcal{C}(\phi)$ due to Ma-Minda \cite{ma minda-1992-a}.
\begin{lem} \label{him-p3-lem-1.6} \cite{ma minda-1992-a}
	Let $f \in \mathcal{S}^{*}(\phi)$. Then $zf'(z)/f(z) \prec zH'(z)/H(z)$ and $f(z)/z \prec H(z)/z$.
\end{lem}
\begin{lem} \label{him-p2-lem-1.13} \cite{ma minda-1992-a}
	Let $f \in \mathcal{C}(\phi)$. Then $zf''(z)/f'(z) \prec zK''(z)/K'(z)$ and $f'(z) \prec K'(z)$.
\end{lem} 
\begin{lem} \label{him-p2-lem-1.14} \cite{ma minda-1992-a}
	Assume $f \in \mathcal{C}(\phi)$ and $|z|=r<1$. Then 
	\begin{equation} \label{him-p2-e-1.15}
	K'(-r) \leq |f'(z)| \leq K'(r).
	\end{equation}
	Equality holds for some $z \neq 0$ if, and only, if f is a rotation of $K$.
\end{lem}
\vspace{4mm}
In \cite{Ravichandran-2004}, Ravichandran has considered the classes $\mathcal{S}^{*}_{c}(\phi)$ and $\mathcal{C}_{c}(\phi)$, the classes of Ma-Minda type starlike functions with respect to the conjugate points and classes of Ma-Minda type convex functions with respect to the conjugate points respectively. A function $f \in \mathcal{\mathcal{S}}$ is in the class $\mathcal{S}^{*}_{c}(\phi)$ if 
$$
\dfrac{2zf'(z)}{f(z)+\overline{f(\bar{z})}} \prec \phi (z)\quad \mbox{for} \quad z\in \mathbb{D}
$$ 
and is in the class $\mathcal{C}_{c}(\phi)$ if 
$$
\dfrac{2(zf'(z))'}{\left(f(z)+\overline{f(\bar{z})}\right)'} \prec \phi (z) \quad \mbox{for} \quad z\in \mathbb{D}.
$$
If $\phi(z)=(1+z)/(1-z)$, then $\mathcal{S}^{*}_{c}(\phi)$ and $\mathcal{C}_{c}(\phi)$ reduce to the classes of standard starlike and convex functions with respect to the conjugate points. 
\begin{lem} \cite{Ravichandran-2004} \label{him-p3-lem-1.14}
	Let $\min _{|z|=r} |\phi (z)|=\phi (-r)$, $\max _{|z|=r}|\phi (z)|=\phi (r)$, $|z|=r$. If $f \in \mathcal{C}_{c}(\phi)$, then 
	\begin{enumerate}
		\item [(i)] $K'(-r) \leq |f'(z)| \leq K'(r)$ \\
		\item[(ii)] $-K(-r) \leq |f(z)| \leq K(r)$ \\
		\item[(iii)] $f(\mathbb{D})\supseteq \{w:|w|\leq -K(-1)\}$.
	\end{enumerate}
	The results are sharp.
\end{lem}

\vspace{2mm}
\noindent Recall that a complex-valued function $f$ in $\mathbb{D}$ is said to be harmonic if it satisfies the Laplace equation $\bigtriangleup f=4f_{z\bar{z}}=0$. Every harmonic mapping $f$ in $\mathbb{D}$ has the unique {\it canonical decomposition} $f=h+\overline{g}$, where $h$ and $g$ are analytic functions with $g(0)=0$. We know that a harmonic mapping $f$ is locally univalent at $z_{0}$ if, and only if, its Jacobian $J_{f}(z)=|h'(z)|^{2}-|g'(z)|^{2} \neq 0$ at $z_{0}$, and is sense-preserving if $J_{f}(z)>0$ in $\mathbb{D}$ {\it i.e.,} the dilation $\omega$ of $f$, given by $\omega(z)=g'(z)/h'(z)$, satisfies $|\omega(z)|<1$ in $\mathbb{D}$. Let $\mathcal{H}$ be the class of normalized harmonic mappings $f=h+\overline{g}$ in $\mathbb{D}$ of the form 
\begin{equation} \label{him-p4-e-1.5}
h(z)=z+ \sum \limits_{n=2}^{\infty} a_{n}z^{n} \quad \mbox{and} \quad  g(z)=\sum \limits_{n=2}^{\infty} b_{n}z^{n}
\end{equation}
and $\mathcal{S}_{\mathcal{H}}$ be the subclass of $\mathcal{H}$ consisting of univalent and sense-preserving harmonic mappings in $\mathbb{D}$. It is known that the class $\mathcal{S}_{\mathcal{H}}$ is normal but not compact. If $g\equiv0$, then the class $\mathcal{S}_{\mathcal{H}}$ reduces to the class $\mathcal{S}$ in $\mathbb{D}$.

\vspace{4mm}

In \cite{Y. Sun}, Sun {\it et.al.} have introduced the class $\mathcal{M}(\alpha, \beta)$ consisting of harmonic mappings $f$ of the form \eqref{him-p4-e-1.5}, with $h'(0) \neq 0$, which satisfy $g'(z)=\alpha z h'(z)$ and $\real \left( 1+zf''(z)/f'(z)\right)> \beta $ in $\mathbb{D}$. The class $\mathcal{M}(1, -1/2)$ with $|\alpha|=1$, has been extended to $\mathcal{M}(\alpha, -1/2)$ by Bharanedhar and Ponnusamy in \cite{bharandar-2014}. Motivated by the class $\mathcal{M}(\alpha, \beta)$, we consider the following new subclasses of $\mathcal{H}$ as follows: 
\begin{defn}
For $\alpha \in \mathbb{C}$ with $|\alpha|\leq 1$, let $\mathcal{HC}(\phi)$ and $\mathcal{HC}_{c}(\phi)$ denote the class of harmonic mappings $f=h+\overline{g}$ in $\mathbb{D}$ of the form \eqref{him-p4-e-1.5}, whose analytic part $h$ belongs $\mathcal{C}(\phi)$ and $ \mathcal{C}_{c}(\phi)$ respectively, with $h'(0)\neq0$, along with the condition $g'(z)=\alpha z h'(z)$.
\end{defn}
 It is easy to see that the dilation of functions belongs to $\mathcal{HC}(\phi)$ and $\mathcal{HC}_{c}(\phi)$ are $\omega(z)=\alpha z $ and $|\omega(z)|<1$. Hence these classes are sense-preserving in $\mathbb{D}$. For a general function $\phi$, it is not easy to obtain the coefficient bounds for the above defined classes. But when $\phi(z)=(1+(1-2 \beta)z)/(1-z)$, the class $\mathcal{HC}(\phi)$ reduces to $\mathcal{M}(\alpha, \beta)$ and the coefficient estimates have been obtained in \cite{Y. Sun}. 
\begin{lem} \cite{Y. Sun} \label{him-p3-lem-1.17}
	Let $f \in \mathcal{M}(\alpha, \beta)$ be of the form \eqref{him-p4-e-1.5}. Then 
	\begin{enumerate}
		\item [(i)] $|a_{n}| \leq \dfrac{1}{n!} \prod \limits _{j=0}^{n} (j-2\beta) \quad (n=2,3 \cdots)$,
		\item[(ii)] $|b_{2}| = \alpha/2$ and $|b_{n}| \leq \dfrac{(n-1)|\alpha|}{n!} \prod \limits _{j=0}^{n} (j-2\beta) \quad (n=3,4, \cdots)$.
	\end{enumerate}
	Moreover, these bounds are sharp with the extremal functions 
	\begin{equation} \label{him-p3-e-1.18}
	f_{\alpha, \beta}(z)= \int \limits _{0}^{z} \dfrac{dt}{(1-\gamma t)^{2-2\beta}} + \overline{\int \limits _{0}^{z} \dfrac{\alpha tdt}{(1-\gamma t)^{2-2\beta}}} \quad (|\gamma|=1; \quad z \in \mathbb{D}).
	\end{equation}
\end{lem}

\begin{lem} \cite{Y. Sun} \label{him-p3-lem-1.19}
	Let $f \in \mathcal{M}(\alpha, \beta)$ with $0 \leq \beta <1$. Then $f$ satisfies the following inequalities 
	\begin{equation} \label{him-p3-e-1.19-a}
	L(r,\alpha,\beta) \leq |f(z)| \leq R(r,\alpha,\beta),
	\end{equation}
	where 
	\[ 
	L(r,\alpha,\beta)=
	\begin{cases}
	\dfrac{(1+|\alpha|)r}{1+r} -|\alpha|log(1+r), \quad &  \beta=0 \\[2mm]
	-|\alpha|r + (1+|\alpha|)log(1+r), \quad &  \beta =1/2\\[2mm]
	\dfrac{-(|\alpha| +2\beta )(1+r)+(1+r)^{2\beta} \left(|\alpha| +2\beta -(2\beta -1)|\alpha|r \right) }{2\beta (2\beta -1) (1+r)}, \quad &  \beta \neq 0, \, 1/2
	\end{cases}
	\]
	and
	\[ 
	R(r,\alpha,\beta)=
	\begin{cases}
	\dfrac{(1+|\alpha|)r}{1-r} +|\alpha|log(1-r), \quad &  \beta=0\\[2mm]
	-|\alpha|r - (1+|\alpha|)log(1-r), \quad &  \beta =1/2\\[2mm]
	\dfrac{(|\alpha| +2\beta )(1-r)-(1-r)^{2\beta} \left(|\alpha| +2\beta +(2\beta -1)|\alpha|r \right) }{2\beta (2\beta -1) (1-r)}, \quad &  \beta \neq 0, \, 1/2.
	\end{cases}
	\]
	All these bounds are sharp, the exremal function is $f_{\alpha,\beta}$ or its rotations, where

	\[ 
	f_{\alpha,\beta}(z)=
	\begin{cases}
	\frac{z}{1-z} +\overline{\alpha (\frac{z}{1-z}+log(1-z))}, \quad &  \beta=0\\[2mm]
	-log(1-z) - \overline{\alpha (z+log(1-z))}, \quad &  \beta =1/2\\[2mm]
	\frac{1-(1-z)^{2\beta -1}}{2\beta -1} + \overline{\frac{\alpha}{2\beta (2\beta -1)}\left[1-(1-z)^{2\beta -1}\left(1+(2\beta -1)z\right)\right]}, \quad &  \beta \neq 0, \, 1/2.
	\end{cases}
	\] 
\end{lem}

In 2018, Bhowmik and Das \cite{bhowmik-2018} proved the following interesting result for subordination classes which will help us to prove our main results.
\begin{lem} \label{him-p2-lem-1.23} \cite{bhowmik-2018}
	Let $f(z)=\sum_{n=0}^{\infty} a_{n}z^{n}$ and $g(z)=\sum_{n=0}^{\infty} b_{n}z^{n}$ be two analytic functions in $\mathbb{D}$ and $g \prec f$, then 
	\begin{equation} \label{him-p2-e-1.24}
	\sum_{n=0}^{\infty} |b_{n}| r^{n} \leq \sum_{n=0}^{\infty} |a_{n}| r^{n}
	\end{equation}
	for $z|=r \leq 1/3.$
\end{lem} 
The present paper is organized as follows. In Section $2$, we state our main results whose proofs will be presented in Section $3$. First we obtain sharp growth estimates for the classes $\mathcal{HC}(\phi)$ and $\mathcal{HC}_{c}(\phi)$ in Theorem $2.1$ and the area of the image $f(\mathbb{D}_{r})$ for $f$ belongs to the aforesaid classes in Theorem $2.3$, where $\mathbb{D}_{r}:=\{z\in \mathbb{D}:|z|=r<1\}$. In Theorem $2.5$ and Theorem $2.12$, we establish the Bohr phenomenon for the classes $\mathcal{HC}(\phi)$ and $\mathcal{HC}_{c}(\phi)$ respectively. As a consequence of Theorem $2.5$, we obtain interesting corollaries and remarks. In Theorem $2.11$, an improved version of Bohr's inequality is obtained for the class $\mathcal{HC}(\phi)$.
\section{Main Results}
\noindent First we prove the sharp growth estimate for the classes $\mathcal{HC}(\phi)$ and $\mathcal{HC}_{c}(\phi)$ which will be useful to prove our main results.
\begin{thm} \label{him-p4-thm-2.1}
Let $f \in \mathcal{HC}(\phi)$ (respectively $\mathcal{HC}_{c}(\phi)$). Then $f$ satisfies the following inequalities 
\begin{equation} \label{him-p4-e-2.2}
L(r,\alpha) \leq |f(z)| \leq R(r,\alpha),
\end{equation}
where
$$
L(r,\alpha)=-K(-r)-|\alpha|\int\limits_{0}^{r} tK'(-t)\,\,dt \quad \mbox{and} \quad R(r,\alpha)=K(r)+|\alpha|\int\limits_{0}^{r}tK'(t)\,\,dt.
$$
The bounds are sharp being the extremal function $f_{\alpha}=h_{\alpha}+\overline{g_{\alpha}}$ with $h_{\alpha}=K$ or its rotations.
\end{thm}
Let $\mathbb{D}_{r}:=\{z\in \mathbb{D}:|z|< r<1\}$.
\begin{thm} \label{him-p4-thm-2.3}
Let $f \in \mathcal{HC}(\phi)$ and $S_{r}$ be the area of the image $f(\mathbb{D}_{r})$. Then the following inequalities hold
\begin{equation} \label{him-p4-e-2.4}
2\pi \int\limits_{0}^{r}t\left(1-|\alpha|^{2}t^{2}\right)(K'(-t))^{2}\,\,dt \leq S_{r} \leq 2\pi \int\limits_{0}^{r}t\left(1-|\alpha|^{2}t^{2}\right)(K'(t))^{2}\,\, dt.
\end{equation}
\end{thm}
\noindent In the following theorem we obtain the Bohr inequality for the class $\mathcal{HC}(\phi)$.
\begin{thm} \label{him-p4-thm-2.5}
Let $f \in \mathcal{HC}(\phi)$ be of the form \eqref{him-p4-e-1.5}. Then the majorant series of $f$ satisfies the following inequality
\begin{equation} \label{him-p4-e-2.6}
|z|+\sum_{n=2}^{\infty} (|a_{n}|+|b_{n}|)|z|^{n} \leq d(f(0),\partial f(\mathbb{D}))
\end{equation}
for $|z|=r\leq \min \{1/3,r_{f}\}$, where $r_{f}$ is the smallest positive root of $R_{\mathcal{C}}(r)=L(1,\alpha)$. Here
$R_{\mathcal{C}}(r)=M_{K}(r)+|\alpha|\int_{0}^{r} tM_{K'}(t)\,\, dt$
and $L(r,\alpha)$ is defined as in Theorem \ref{him-p4-thm-2.1}. 
\end{thm}
\begin{rem} \label{him-p4-rem-2.1}
\begin{enumerate}
\item [(i)] $\left(\mbox{Bohr radius for $\mathcal{HC}(\phi)$ when $\phi$ has positive coefficients}\right)$ \\
Let $\phi(z)=1+\sum_{n=1}^{\infty} B_{n}z^{n}$. It is easy to see that if we impose an additional condition on $\phi$ that the coefficients $B_{n}$'s are positive, then the majorant series is $M_{\phi}(r)=\phi(r)$. From the definition of $H$ given in \eqref{him-p3-e-1.6}, we obtain 
\begin{equation} \label{him-p2-e-3.8}
H(z)= z \exp \left(\int\limits_{0}^{z} \dfrac{\phi(t)-1}{t}\,\, dt\right)=z \exp \left(\sum \limits _{n=1}^{\infty} \frac{B_{n}}{n}z^{n}\right).
\end{equation}
From \eqref{him-p2-e-3.8}, it is easy to see that all the coefficients of $H$ are positive and hence it follows that $M_{H}(r)=H(r)$. Thus, the relation $zK'(z)=H(z)$ yields that $K(z)$ has all positive coefficients and hence $M_{K}(r)=K(r)$. Therefore, we obtain 
\begin{equation} \label{him-p2-e-3.8-a}
 M_{K'}(r)=K'(r)\quad \mbox{and} \quad R_{\mathcal{C}}(r)=K(r)+|\alpha| \int\limits_{0}^{r} t K'(t) \,\, dt=:R(r).
\end{equation}  
Then each $f \in \mathcal{HC}(\phi)$ satisfies the inequality \eqref{him-p4-e-2.6} for $|z|\leq \min \{1/3,r_{f}\}$, where $r_{f}$ is the root of the equation $R(r)=L(1,\alpha)$. In particular, when $r_{f} \leq 1/3$, the radius $r_{f}$ is the best possible for the function $f_{\alpha}=h_{\alpha}+\overline{g_{\alpha}}$ with $h_{\alpha}=K$ or its rotations, since $K\in \mathcal{C}(\phi)$ with real positive coefficients. Indeed, for $|z|=r_{f}$, 
$$
M_{f_{\alpha}}(r_{f})=R(r_{f})=L(1,\alpha)=d(f_{\alpha}(0),\partial f_{\alpha}(\mathbb{D})),
$$
which shows that $r_{f}$ is the best possible.
\item [(ii)] (Bohr phenomenon for corresponding class $\mathcal{HC}(\psi)$ associated with non-Ma-Minda function) 
Let $\psi$ be the corresponding non-Ma-Minda function of $\phi$. Since $\psi$ is actually obtained from $\phi$ by a rotation $z$ by $-z$, the image of the unit disk $\mathbb{D}$ under the functions $\psi$ and $\phi$ are identical. Therefore, we conclude that $\mathcal{HC}(\psi)=\mathcal{HC}(\phi)$ and the Bohr radius for the class $\mathcal{HC}(\psi)$ is same as that of $\mathcal{HC}(\phi)$.

\item [(iii)] Observe that $\mathcal{HC}(\phi)$ reduces to $\mathcal{C}(\phi)$, if the co-analytic part of $f$ vanishes identically in $\mathbb{D}$. Thus, as a consequence of Theorem \ref{him-p4-thm-2.5}, we can obtain the radius $r_{f}$ so that Bohr inequality \eqref{him-p4-e-1.3} holds for $|z|=r\leq r_{f}$, which has been extensively studied by Allu and Halder \cite{Himadri-Vasu-P2}. For particular values of $\phi$, Allu and Halder \cite{Himadri-Vasu-P2} have obtained the Bohr radius $r_{f}$.
\end{enumerate}
\end{rem}
As a consequence of Theorem \ref{him-p4-thm-2.1} and Remark \ref{him-p4-rem-2.1}, for $\phi(z)=1+4z/3 +2z^{2}/3$, we obtain the following corollary.
\begin{cor} \label{him-p4-cor-2.9}
Let $0.53143<|\alpha|<1$. Then every function $f \in \mathcal{HC}(\phi)$ of the form \eqref{him-p4-e-1.5} with $\phi(z)=1+4z/3 +2z^{2}/3$ satisfies the inequality \eqref{him-p4-e-2.6} for $|z|=r \leq r_{f}\leq 1/3$. The constant $r_{f}$ is the best possible.
\end{cor}
In view of Theorem \ref{him-p4-thm-2.1}, Remark \ref{him-p4-rem-2.1} and going by the similar lines of arguments as in the proof of Corollary \ref{him-p4-cor-2.9}, we can establish the Bohr phenomenon for other particular choices of $\phi$. For $\phi(z)=(1+(1-2 \beta)z)/(1-z)$ with $-1/2 \leq \beta <1$, the class $\mathcal{HC}(\phi)$ reduces to $\mathcal{M}(\alpha, \beta)$. Using Theorem \ref{him-p4-thm-2.1}, one can obtain the radius $r_{f}$ for the class $\mathcal{M}(\alpha, \beta)$ so that the inequality \eqref{him-p4-e-2.6} holds. However, in that case when $r_{f}>1/3$, we cannot ensure that $r_{f}$ is the best possible. Using the sharp coefficient estimates for $\mathcal{M}(\alpha, \beta)$ (see \cite{Y. Sun}), we obtain the sharp radius $r_{f}$ in the following corollary.

\begin{cor} \label{him-p4-cor-2.10}
	Let $f \in \mathcal{M}(\alpha,\beta)$ be of the form \eqref{him-p4-e-1.5} with $|\alpha|\leq 1, 0 \leq \beta <1$. Then the inequality \eqref{him-p4-e-2.6} satisfies for $|z|=r \leq r_{f}$, where $r_{f}$ is the smallest root of $D_{1}(r):=R(r,\alpha,\beta)-L(1,\alpha,\beta)=0$. The radius $r_{f}$ is sharp. 
\end{cor}

The roots $r_{f}$ of $D_{1}(r)=0$ for different values of $\alpha$ and $\beta$ have been shown in Tables $1,2,3$ and Figures $1,2$.
\vspace{3mm}

		

\begin{figure}[!htb]
	\begin{center}
		\includegraphics[width=0.50\linewidth]{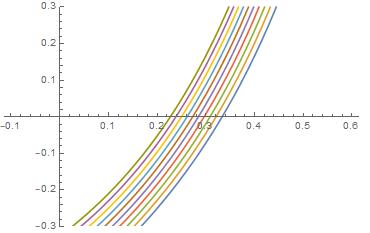}
		\;\;\;
		\includegraphics[width=0.40\linewidth]{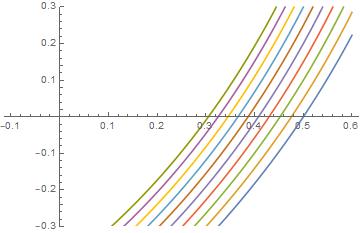}
	\end{center}
	\caption{The graphs of $D_{1}(r)$ respectively for $\beta = 0$ and $\beta = 1/2$ when $\alpha=0.0,\,  0.1, 0.2,\, 0.3,$ $ 0.4,\, 0.5,\, 0.6,\, 0.7,\, 0.8,\, 0.9$.}
	\label{figure-2}
\end{figure} 

\begin{table}[ht]
	
	\begin{tabular}{|l|l|l|l|l|l|l|l|l|l|l|} 
		\hline
		$\alpha$ & $0.0$ & $0.1$ & $0.2$ & $0.3$ & $0.4$ & $0.5$ & $0.6$& $0.7$ & $0.8$ & $0.9$  \\
		\hline
		$\beta$ & $0.0$ & $0.0$ & $0.0$& $0.0$ & $0.0$ & $0.0$ & $0.0$ & $0.0$ & $0.0$ & $0.0$\\
		\hline
		$r_{f}$ & $0.333$ & $0.321$ & $0.308$& $0.296$ & $0.284$ & $0.273$ & $0.261$ & $0.25$ & $0.238$ & $0.227$\\
		\hline
	\end{tabular}
	\vspace{1mm}
	
	\caption{The roots $r_{f}$ of $D_{1}(r)=0$ for different values of $\alpha$ when $\beta=0$.}
\end{table}

\begin{table}[ht]
	
	\begin{tabular}{|l|l|l|l|l|l|l|l|l|l|l|} 
		\hline
		$\alpha$ & $0.0$ & $0.1$ & $0.2$ & $0.3$ & $0.4$ & $0.5$ & $0.6$& $0.7$ & $0.8$ & $0.9$  \\
		\hline
		$\beta$ & $0.5$ & $0.5$ & $0.5$& $0.5$ & $0.5$ & $0.5$ & $0.5$ &   $0.5$ & $0.5$ & $0.5$\\
		\hline
		$r_{f}$ & $0.5$ & $0.476$ & $0.452$& $0.43$ & $0.408$ & $0.387$ & $0.366$ & $0.345$ & $0.321$ & $0.305$\\
		\hline
	\end{tabular}
	\vspace{1mm}
	
	\caption{The roots $r_{f}$ of $D_{1}(r)=0$ for different values of $\alpha$ when $\beta=1/2$.}
\end{table}


\begin{table}[ht]
	
	\begin{tabular}{|l|l|l|l|l|l|l|l|l|l|l|} 
		\hline
		$\alpha$ & $0.0$ & $0.1$ & $0.2$ & $0.3$ & $0.4$ & $0.5$ & $0.6$& $0.7$ & $0.8$ & $0.9$  \\
		\hline
		$\beta$ & $0.9$ & $0.9$ & $0.9$& $0.9$ & $0.9$ & $0.9$ & $0.9$ & $0.9$ & $0.9$ & $0.9$\\
		\hline
		$r_{f}$ & $0.815$ & $0.757$ & $0.705$& $0.656$ & $0.61$ & $0.568$ & $0.527$ & $0.488$ & $0.451$ & $0.415$\\
		\hline
	\end{tabular}
	\vspace{1mm}
	
	\caption{The roots $r_{f}$ of $D_{1}(r)=0$ for different values of $\alpha$ when $\beta=0.9$.}
\end{table}

\begin{figure}[!htb]
	\begin{center}
		\includegraphics[width=0.450\linewidth]{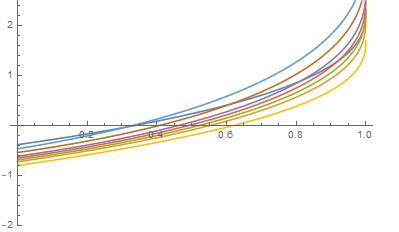}
		\;\;\;\;\;
		\includegraphics[width=0.450\linewidth]{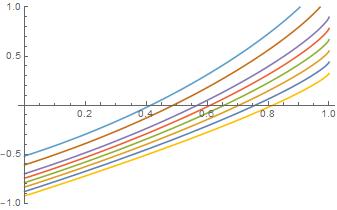}
	\end{center}
	\caption{The graphs of $D_{1}(r)$ respectively for $\beta=0.7$ and $\beta=0.9$ when $\alpha=0.0,\,  0.1, 0.2,\,$  $0.3, 0.4,\, 0.5,\, 0.6,\, 0.7,\, 0.8,\, 0.9$.}
	\label{figure-3}
\end{figure}


\vspace{3mm}

\noindent In the next Theorem, we establish an improved version of Bohr inequality for the class $\mathcal{HC}(\phi)$ and show that after adding area quantity $S_{r}/2\pi$ with the Majorant series of $f \in \mathcal{HC}(\phi)$, the sum is still less than $d(f(0),\partial f(\mathbb{D}))$ for some radius $r \leq R^{'}_{f}<1$.
\begin{thm} \label{him-p4-thm-2.11}
	Let $f \in \mathcal{HC}(\phi)$ be of the form \eqref{him-p4-e-1.5}. Let $S_{r}$ be the area of the image $f(\mathbb{D}_{r})$. Then the following inequality 
	$$
	M_{f}(r)+ \frac{S_{r}}{2\pi} \leq d(f(0),\partial f(\mathbb{D}))
	$$
	holds for $|z|=r\leq R^{'}_{f}:= \min \{1/3,r'_{f}\}$, where $r'_{f}$ is the smallest positive root of $R'_{f}(r)=L(1,\alpha)$. Here 
	$R'_{f}(r)=M_{K}(r)+|\alpha|\int_{0}^{r} tM_{K'}(t)\,\, dt + \int_{0}^{r}t\left(1-|\alpha|^{2}t^{2}\right)(K'(t))^{2}\,\, dt$ and $L(r,\alpha)$ is defined as in Theorem \ref{him-p4-thm-2.1}.
\end{thm}

\noindent Next, we establish the Bohr phenomenon for the class $\mathcal{HC}_{c}(\phi)$.
\begin{thm} \label{him-p4-thm-2.7}
Let $f \in \mathcal{HC}_{c}(\phi)$ be of the form \eqref{him-p4-e-1.5}. Then the majorant series of $f$ satisfies the following inequality
\begin{equation} \label{him-p4-e-2.8}
|z|+\sum_{n=2}^{\infty} (|a_{n}|+|b_{n}|)|z|^{n} \leq d(f(0),\partial f(\mathbb{D}))
\end{equation}
for $|z|=r\leq \min \{1/3,r_{f}\}$, where $r_{f}$ is the smallest positive root of $R_{\mathcal{C}_{c}}(r)=L(1,\alpha)$. Here
$R_{\mathcal{C}_{c}}(r)=T(r)+|\alpha|\int_{0}^{r} tT_{c}(t)\,\, dt$ with 
\begin{equation} \label{him-p4-e-2.8-a}
T_{c}(r):=\frac{1}{r}\int \limits _{0}^{r}M_{K'}(t) M_{\phi}(t) \, dt,\quad T(r)=\int\limits_{0}^{r}T_{c}(t)\,\, dt 
\end{equation}
 and $L(r,\alpha)$ is defined as in Theorem \ref{him-p4-thm-2.1}. 
\end{thm}
\section{Proof of the main results}
\begin{pf} [{\bf Proof of Theorem   \ref{him-p4-thm-2.1}}]
Let $f=h+\overline{g} \in \mathcal{HC}(\phi)$ (respectively $\mathcal{HC}_{c}(\phi)$). Then from definition, we have $h \in \mathcal{C}(\phi)$ (respectively $\mathcal{C}_{c}(\phi)$). In view of Lemmas \ref{him-p2-lem-1.14} and \ref{him-p3-lem-1.14}, we obtain 
\begin{equation} \label{him-p4-e-3.1}
K'(-r) \leq |h'(z)| \leq K'(r) \quad \mbox{for} \quad |z|=r.
\end{equation}
Let $\gamma$ be the linear segment joining $0$ to $z$ in $\mathbb{D}$. Then we have 
\begin{align} \label{him-p4-e-3.2}
|f(z)| &=\left|\int_{\gamma}\frac{\partial f}{\partial \xi}\,\, d\xi +\frac{\partial f}{\partial \bar{\xi}}\,\, d\bar{\xi} \right| \\[2mm] \nonumber
& \leq \int_{\gamma} \left(|h'(\xi)|+|g'(\xi)|\right)\,\, |d\xi|\\[2mm] \nonumber
&=\int_{\gamma} \left(1+|\alpha||\xi|\right)|h'(\xi)|\,\, |d\xi|.
\end{align}
Hence by using \eqref{him-p4-e-3.1} and \eqref{him-p4-e-3.2}, we obtain 
\begin{equation} \label{him-p4-e-3.3}
|f(z)|\leq \int \limits_{0}^{r} \left(1+|\alpha|t\right)K'(t) \,\, dt=K(r)+|\alpha|\int\limits_{0}^{r}tK'(t)\,\,dt= R(r,\alpha).
\end{equation}
Let $\Gamma$ be the preimage of the line segment joinig $0$ to $f(z)$ under the function $f$. Then we have
\begin{align} \label{him-p4-e-3.4}
|f(z)|&=\left|\int_{\Gamma}\frac{\partial f}{\partial \xi}\,\, d\xi +\frac{\partial f}{\partial \bar{\xi}}\,\, d\bar{\xi} \right|\\[2mm] \nonumber 
& \geq \int_{\Gamma} \left(|h'(\xi)|-|g'(\xi)|\right)\,\, |d\xi| \\[2mm] \nonumber &=\int_{\Gamma} \left(1-|\alpha||\xi|\right)|h'(\xi)|\,\, |d\xi|.
\end{align}
In view of \eqref{him-p4-e-3.1} and \eqref{him-p4-e-3.4}, we obtain 
\begin{equation} \label{him-p4-e-3.5}
|f(z)|\geq \int \limits_{0}^{r} \left(1-|\alpha|t\right)K'(-t) \,\, dt=-K(-r)-|\alpha|\int\limits_{0}^{r}tK'(-t)\,\,dt=L(r,\alpha).
\end{equation}
From \eqref{him-p4-e-3.3} and \eqref{him-p4-e-3.5}, we have 
\begin{equation} \label{him-p4-e-3.6}
 L(r,\alpha) \leq |f(z)| \leq R(r,\alpha).
\end{equation}
To show the sharpness we consider the function $f_{\alpha}=h_{\alpha}+\overline{g_{\alpha}}$ with $h_{\alpha}=K$ or its rotations. It is easy to see that $h_{\alpha}=K \in \mathcal{C}(\phi)$ and satisfies $g'_{\alpha}(z)=\alpha z h'_{\alpha}(z)$, which shows that $f_{\alpha}\in \mathcal{HC}(\phi)$. The equality holds on both sides of \eqref{him-p4-e-3.1} for suitable rotations of $K$. For $0<\alpha<1$, it is easy to see that $f_{\alpha}(r)=R(r,\alpha)$ and $f_{\alpha}(-r)=-L(r,\alpha)$. Hence $|f_{\alpha}(r)|=R(r,\alpha)$ and $|f_{\alpha}(-r)|=L(r,\alpha)$. This completes the proof.
\end{pf}
\begin{pf} [{\bf Proof of Theorem   \ref{him-p4-thm-2.3}}]
Let $f=h+\overline{g} \in \mathcal{HC}(\phi)$ and $z=x+iy$. Then the area of image of $\mathbb{D}_{r}$ under a harmonic function $f$ is given by 
\begin{equation} \label{him-p4-e-3.7}
S_{r}=\iint_{\mathbb{D}_{r}} \left(|h'(z)|^{2}-|g'(z)|^{2}\right)\,\, dx\, dy= \iint_{\mathbb{D}_{r}} \left(1-|\alpha|^{2}|z|^{2}\right)|h'(z)|^{2} \,\, dx\, dy.
\end{equation}	
Since $h\in \mathcal{C}(\phi)$, in view of \eqref{him-p4-e-3.1} and \eqref{him-p4-e-3.7}, we obtain 
$$
 \int\limits_{0}^{r} \int \limits_{\theta=0}^{2\pi}t\left(1-|\alpha|^{2}t^{2}\right)(K'(-t))^{2}\,\,d \theta \,\,dt \leq S_{r} \leq  \int\limits_{0}^{r}\int \limits_{\theta=0}^{2\pi}t\left(1-|\alpha|^{2}t^{2}\right)(K'(t))^{2}\,\, d\theta \,\, dt,
$$
which is equivalent to
$$
2\pi \int\limits_{0}^{r}t\left(1-|\alpha|^{2}t^{2}\right)(K'(-t))^{2}\,\,dt \leq S_{r} \leq 2\pi \int\limits_{0}^{r}t\left(1-|\alpha|^{2}t^{2}\right)(K'(t))^{2}\,\, dt.
$$
This completes the proof.	
\end{pf}
\begin{pf} [{\bf Proof of Theorem   \ref{him-p4-thm-2.5}}]
Let $f=h+\overline{g} \in \mathcal{HC}(\phi)$. Since $h \in \mathcal{C}(\phi)$, from Lemma \ref{him-p2-lem-1.13}, we have 
\begin{equation} \label{him-p4-e-3.8}
h' \prec K'.
\end{equation}
Let $K(z)=z+\sum \limits_{n=2}^{\infty} k_{n}z^{n}$.
In view of Lemma \ref{him-p2-lem-1.23} and \eqref{him-p4-e-3.8}, we obtain 
\begin{equation}\label{him-p4-e-3.9}
1+\sum\limits_{n=2}^{\infty} n|a_{n}|r^{n-1}=M_{h'}(r) \leq M_{K'}(r)=1+\sum\limits_{n=2}^{\infty} n|k_{n}|r^{n-1}
\end{equation}
for $|z|=r\leq 1/3$. Integrating \eqref{him-p4-e-3.9} with respect to $r$ from $0$ to $r$, we obtain 
\begin{equation} \label{him-p4-e-3.10}
M_{h}(r)=r+\sum\limits_{n=2}^{\infty} |a_{n}|r^{n} \leq r+\sum\limits_{n=2}^{\infty} |k_{n}|r^{n}=M_{K}(r) \quad \mbox{for} \quad r\leq 1/3.
\end{equation}
From the definition of $\mathcal{HC}(\phi)$, we have $g'(z)=\alpha z h'(z)$. This relation along with \eqref{him-p4-e-3.9} gives
\begin{equation} \label{him-p4-e-3.11}
\sum\limits_{n=2}^{\infty} n|b_{n}|r^{n-1}=M_{g'}(r)=|\alpha|r M_{h'}(r) \leq |\alpha|r M_{K'}(r)\quad \mbox{for} \quad r\leq 1/3.
\end{equation} 
By integrating \eqref{him-p4-e-3.11} with respect to $r$ from $0$ to $r$, we obtain 
\begin{equation} \label{him-p4-e-3.12}
M_{g}(r)=\sum\limits_{n=2}^{\infty} |b_{n}|r^{n} \leq |\alpha| \int\limits_{0}^{r} t M_{K'}(t) \,\, dt\quad \mbox{for} \quad r\leq 1/3.
\end{equation}
Therefore, for $|z|=r\leq 1/3$, the inequalities \eqref{him-p4-e-3.10} and \eqref{him-p4-e-3.12} yeild that 
\begin{equation} \label{him-p4-e-3.13}
M_{f}(r)=|z|+\sum_{n=2}^{\infty} (|a_{n}|+|b_{n}|)r^{n} \leq M_{K}(r)+ |\alpha| \int\limits_{0}^{r} t M_{K'}(t) \,\, dt=R_{\mathcal{C}}(r).
\end{equation}
From the inequality \eqref{him-p4-e-2.2}, it is evident that the Euclidean distance between $f(0)$ and the boundary of $f(\mathbb{D})$ is given by 
\begin{equation} \label{him-p4-e-3.14}
d(f(0), \partial f(\mathbb{D}))= \liminf \limits_{|z|\rightarrow 1} |f(z)-f(0)| \geq L(1,\alpha).
\end{equation}
 We note that $R_{\mathcal{C}}(r) \leq L(1)$ whenever $r \leq r_{f}$, where $r_{f}$ is the smallest positive root of $R_{\mathcal{C}}(r)=L(1,\alpha)$ in $(0,1)$. Let $H_{1}(r)=R_{\mathcal{C}}(r)-L(1,\alpha)$ then $H_{1}(r)$ is a continuous function in $[0,1]$. Since $M_{K}(r) \geq K(r)>-K(-r)$, it follows that 
\begin{align} \label{him-p4-e-3.15}
H_{1}(1)&=R_{\mathcal{C}}(1)-L(1,\alpha)\\ \nonumber
&=  M_{K}(1)+K(-1)+|\alpha|\int\limits_{0}^{r} t \left(M_{K'}(t)+K'(t)\right) \,\, dt \\ \nonumber
&\geq  K(1)+K(-1)+|\alpha|\int\limits_{0}^{r} t \left(M_{K'}(t)+K'(t)\right) \,\, dt>0.
\end{align}
On the other hand,
\begin{equation} \label{him-p4-e-3.15-a}
H_{1}(0)=-L(1,\alpha)=-K(-1)(1-|\alpha|)+|\alpha|\int \limits_{0}^{1}-K(-t)\,\, dt<0.
\end{equation}
Therefore, $H_{1}$ has a root in $(0,1)$. Let $r_{f}$ be the smallest root of $H_{1}$ in $(0,1)$. Then $R_{\mathcal{C}}(r)\leq L(1,\alpha)$ for $r\leq r_{f}$.	Now by combining the inequalities \eqref{him-p4-e-3.13} and \eqref{him-p4-e-3.14} with the fact that $R_{\mathcal{C}}(r)\leq L(1,\alpha)$ for $r\leq r_{f}$, we obtain 
$$
|z|+\sum_{n=2}^{\infty} (|a_{n}|+|b_{n}|)r^{n} \leq d(f(0), \partial f(\mathbb{D}))
$$
for $|z|=r\leq \min \{1/3, r_{f}\}$. This completes the proof.
\end{pf}
\begin{pf} [{\bf Proof of Corollary   \ref{him-p4-cor-2.9}}]
Let $f \in \mathcal{HC}(\phi)$ be of the form \eqref{him-p4-e-1.5} with $\phi(z)=1+4z/3 +2z^{2}/3$. From the relation \eqref{him-p3-e-1.6}, we have 
\begin{equation}
K'(r)=\exp \left(\frac{4}{3}r +\frac{r^{2}}{3}\right),\quad K(r)=\int\limits_{0}^{r}\exp \left(\frac{4}{3}t +\frac{t^{2}}{3}\right)\,\, dt. 
\end{equation}
A simple computation shows that $K(1/3)\approx 0.425549$ and $K(-1)\approx -0.598691$. Further, a simple computation using Mathematica shows that 
\begin{equation} \label{him-p4-e-3.16}
 \int\limits_{0}^{\frac{1}{3}} tK'(t)\,\,dt=0.0766\quad \mbox{and}\quad \int\limits_{0}^{1} (-t)K'(-t) \,\,dt=-0.249202.
\end{equation}
Since all the coefficients of $\phi$ are positive, in view of Remark \ref{him-p4-rem-2.1} (i), we obtain $H_{1}(r)=R(r)-L(1,\alpha)$. A simple computation using \eqref{him-p4-e-3.16} shows that 
$$
H_{1}(0)=-L(1,\alpha)=-0.598691+|\alpha|(0.249202)<0.
$$
On the other hand,
$$
H_{1}(1/3)=R(1/3)-L(1,\alpha)=0.425549+|\alpha|(0.0766)-0.598691+|\alpha|(0.249202)>0
$$
provided $0.5314332 <|\alpha|<1$. Since $H_{1}(0)<0$ and $H_{1}(1/3)>0$, $H_{1}(r)$ has a root in $(0,1/3)$ and choose $r_{f}$ to be the smallest root in $(0,1/3)$. By Remark \ref{him-p4-rem-2.1}, the radius $r_{f}$ is the best possible.	
\end{pf}
\begin{pf} [{\bf Proof of Corollary   \ref{him-p4-cor-2.10}}]
From Lemma \ref{him-p3-lem-1.19}, it is evident that the Euclidean distance between $f(0)$ and the boundary of $f(\mathbb{D})$ is 
\begin{equation} \label{him-p3-e-3.14}
d(f(0), \partial f(\mathbb{D}))= \liminf \limits_{|z|\rightarrow 1} |f(z)-f(0)|\geq L(1,\alpha,\beta).
\end{equation}
We note that $r_{f}$ is the root of the equation $R(r,\alpha,\beta)=L(1,\alpha,\beta)$ in $(0,1)$. The existance of the root is ensured by the relation $R(1,\alpha,\beta) >L(1,\alpha,\beta)$ in view of \eqref{him-p3-e-1.19-a}. For $0<r\leq r_{f}$, it is easy to see that $R(r,\alpha,\beta)\leq L(1,\alpha,\beta)$. In view of Lemma \ref{him-p3-lem-1.17} and \eqref{him-p3-e-3.14}, for $|z|=r\leq r_{f}$, we obtain 
\begin{align*}
|z|+\sum\limits_{n=2}^{\infty} (|a_{n}|+|b_{n}|)|z|^{n} 
&\leq r_{f} + (|a_{2}|+|b_{2}|)r_{f}^{2}+\sum\limits_{n=3}^{\infty} (|a_{n}|+|b_{n}|)r_{f}^{n}\\
& =R(r_{f},\alpha,\beta) \leq  L(1,\alpha,\beta) \leq d(f(0), \partial f(\mathbb{D})).
\end{align*}
To show the sharpness of the radius $r_{f}$, we consider the function $f=f_{\alpha,\beta}$, which is defined in Lemma \ref{him-p3-lem-1.19}. It is easy to see that $f_{\alpha,\beta}$ belongs to $\mathcal{M}(\alpha,\beta)$. Since the left side of the growth inequality in Lemma \ref{him-p3-lem-1.19} holds for $f=f_{\alpha,\beta}$ or its rotations, we have $d(f(0), \partial f(\mathbb{D}))=L(1,\alpha,\beta)$. Therefore, the function $f=f_{\alpha,\beta}$ for $|z|=r_{f}$ gives
\begin{align*}
|z|+\sum\limits_{n=2}^{\infty} (|a_{n}|+|b_{n}|)|z|^{n} 
&= r_{f} + (|a_{2}|+|b_{2}|)r_{f}^{2}+\sum\limits_{n=3}^{\infty} (|a_{n}|+|b_{n}|)r_{f}^{n}\\
& =R(r_{f},\alpha,\beta) = L(1,\alpha,\beta) =d(f(0), \partial f(\mathbb{D})),
\end{align*}
which shows that the radius $r_{f}$ is the best possible. This completes the proof.	
\end{pf}

\begin{pf} [{\bf Proof of Theorem   \ref{him-p4-thm-2.11}}]
Let $f \in \mathcal{HC}(\phi)$ be of the form \eqref{him-p4-e-1.5}.	Then, from the right hand inequality in \eqref{him-p4-e-2.4} and \eqref{him-p4-e-3.13}, we obtain 
\begin{align} \label{him-p4-e-3.20-a}
M_{f}(r)+ \frac{S_{r}}{2\pi} &\leq M_{K}(r)+|\alpha|\int\limits_{0}^{r} tM_{K'}(t)\,\, dt + \int\limits_{0}^{r}t\left(1-|\alpha|^{2}t^{2}\right)(K'(t))^{2}\,\, dt \\ \nonumber
&= R_{\mathcal{C}}(r) +\int\limits_{0}^{r}t\left(1-|\alpha|^{2}t^{2}\right)(K'(t))^{2}\,\, dt =R'_{f}(r)
\end{align}
for $r\leq 1/3$. Let $H_{2}(r)=R'_{f}(r)-L(1,\alpha)$, then $H_{2}(r)$ is a continuous function in $[0,1]$. The inequality \eqref{him-p4-e-3.15-a} yields that $H_{2}(0)=-L(1,\alpha)<0$. From \eqref{him-p4-e-3.15}, we have 
\begin{equation} \label{him-p4-e-3.20-b}
R_{\mathcal{C}}(1)-L(1,\alpha)>0.
\end{equation}
 Observe that, for $|\alpha|<1$, the quantity  $t\left(1-|\alpha|^{2}t^{2}\right)(K'(t))^{2}\geq 0$ and hence 
\begin{equation} \label{him-p4-e-3.20-c}
\int \limits_{0}^{r}t\left(1-|\alpha|^{2}t^{2}\right)(K'(t))^{2}\,\, dt>0.
\end{equation}
In view of the inequalities \eqref{him-p4-e-3.20-a} and \eqref{him-p4-e-3.20-b}, we obtain 
$$
H_{2}(1)=R_{\mathcal{C}}(1)-L(1,\alpha) + \int\limits_{0}^{1}t\left(1-|\alpha|^{2}t^{2}\right)(K'(t))^{2}\,\, dt>0.
$$
Since $H_{2}(0)<0$ and $H_{2}(1)>0$, $H_{2}$ has a root in $(0,1)$ and choose $r'_{f}$ to be the smallest root in $(0,1)$. Theefore, $R'_{f}(r)\leq L(1,\alpha)$ for $r \leq r'_{f}$. Therefore, from the inequality \eqref{him-p4-e-3.14} and \eqref{him-p4-e-3.20-a}, we obtain 
$$
M_{f}(r)+ \frac{S_{r}}{2\pi}\leq d(f(0), \partial f(\mathbb{D}))
$$
for $r \leq \min \{1/3, r'_{f}\}$.		
\end{pf}
\begin{pf} [{\bf Proof of Theorem   \ref{him-p4-thm-2.7}}]
Let $f=h+\overline{g} \in \mathcal{HC}_{c}(\phi)$. Then $h \in \mathcal{C}_{c}(\phi)$. Let $g_{c}(z):=(h(z)+\overline{h(\bar{z})})/2$. Since $\phi$ is starlike and symmetric with respect to real axis, $g_{c} \in \mathcal{C}(\phi)$. From the definition of $\mathcal{C}_{c}(\phi)$, we have 
\begin{equation} \label{him-p3-e-3.13-f}
\left(zh'(z)\right)'=g_{c}'(z)\phi (\omega (z)),
\end{equation}
where $\omega:\mathbb{D}\rightarrow \mathbb{D}$ is analytic with $\omega(0)=0$. A simplication of \eqref{him-p3-e-3.13-f} gives
\begin{equation} \label{him-p3-e-3.13-g}
h'(z)= \frac{1}{z} \int \limits _{0}^{z} g_{c}'(\xi) \phi(\omega(\xi))\, \, d\xi.
\end{equation}
Since $g_{c}\in \mathcal{C}(\phi)$, from Lemma \ref{him-p2-lem-1.13}, we have $g_{c}' \prec K'$ and hence by Lemma \ref{him-p2-lem-1.23}, we obtain 
\begin{equation} \label{him-p3-e-3.13-h}
M_{g_{c}'}(r)\leq M_{K'}(r) \quad \mbox{for } \quad r \leq 1/3.
\end{equation}
Since $\phi \circ \omega \prec \phi$, an application of Lemma \ref{him-p2-lem-1.23} shows that 
\begin{equation} \label{him-p3-e-3.2-b}
M_{\phi \circ \omega} (r) \leq M_{\phi} (r) \quad \mbox{for}\quad |z|=r\leq 1/3.
\end{equation}
In view of \cite[Lemma 2.1]{Himadri-Vasu-P3} and by using \eqref{him-p3-e-3.13-g}, \eqref{him-p3-e-3.13-h} and \eqref{him-p3-e-3.2-b}, we obtain 
\begin{align} \label{him-p4-e-3.19}
M_{h'}(r) &\leq \frac{1}{r}\int \limits _{0}^{r}M_{g_{c}'}(t) M_{\phi\circ \omega}(t) \, dt \\ \nonumber 
&\leq \frac{1}{r}\int \limits _{0}^{r}M_{K'}(t) M_{\phi}(t) \, dt\\ \nonumber 
&=:T_{c}(r)
\end{align}
for $r \leq 1/3$. By integrating \eqref{him-p4-e-3.19} with respect to $r$ from $0$ to $r$, we obtain 
\begin{equation} \label{him-p4-e-3.19-a}
M_{h}(r) \leq \int \limits _{0}^{r} T_{c}(t)\,\, dt=:T(r)\quad \mbox{for} \quad r\leq 1/3. 
\end{equation}	
From the definition of $\mathcal{HC}_{c}(\phi)$, we have $g'(z)=\alpha z h'(z)$. This relation along with the inequality \eqref{him-p4-e-3.19} asserts that 
\begin{equation} \label{him-p4-e-3.20}
\sum\limits_{n=2}^{\infty} n|b_{n}|r^{n-1}=M_{g'}(r)=|\alpha|r M_{h'}(r) \leq |\alpha|r T_{c}(r)\quad \mbox{for} \quad r\leq 1/3.
\end{equation}
Integrating \eqref{him-p4-e-3.20} with respect to $r$ from $0$ to $r$, we obtain
\begin{equation} \label{him-p4-e-3.21}
\sum\limits_{n=2}^{\infty} |b_{n}|r^{n}=M_{g}(r) \leq |\alpha|\int\limits_{0}^{r}tT_{c}(t)\quad \mbox{for} \quad r\leq 1/3.
\end{equation}
Therefore, from the inequalities \eqref{him-p4-e-3.19-a} and \eqref{him-p4-e-3.21}, it follows that 
\begin{equation} \label{him-p4-e-3.22}
M_{f}(r)=|z|+\sum_{n=2}^{\infty} (|a_{n}|+|b_{n}|)r^{n} \leq T(r)+ |\alpha| \int\limits_{0}^{r} t T_{c}(t) \,\, dt=R_{\mathcal{C}_{c}}(r)
\end{equation}
for $|z|=r\leq 1/3$. From the inequality \eqref{him-p4-e-2.2}, it is evident that the Euclidean distance between $f(0)$ and the boundary of $f(\mathbb{D})$ is given by 
\begin{equation} \label{him-p4-e-3.23}
d(f(0), \partial f(\mathbb{D}))= \liminf \limits_{|z|\rightarrow 1} |f(z)-f(0)| \geq L(1,\alpha).
\end{equation}
We note that $R_{\mathcal{C}_{c}}(r) \leq L(1,\alpha)$ whenever $r \leq r_{f}$, where $r_{f}$ is the smallest positive root of $R_{\mathcal{C}_{c}}(r)=L(1,\alpha)$ in $(0,1)$. Let $H_{2}(r)=R_{\mathcal{C}_{c}}(r)-L(1,\alpha)$, then $H_{2}(r)$ is a continuous function in $[0,1]$. Clearly, 
\begin{align} \label{him-p4-e-3.25}
H_{2}(1)&=R_{\mathcal{C}_{c}}(1)-L(1,\alpha)\\[2mm] \nonumber
&= T(1)+K(-1)+|\alpha|\int\limits_{0}^{1} t \left(T'(t)+K'(-t)\right) \,\, dt. 
\end{align}
A simple observation shows that 
\begin{equation}  \label{him-p4-e-3.26}
M_{K'}(r)\geq K'(r)\geq K'(-r),\quad M_{\phi}(r)\geq \phi(r)
\end{equation}
and
\begin{equation} \label{him-p4-e-3.26-a}
 K'(r)+rK''(r)=K'(r) \phi(r)
\end{equation}
Therefore, by using \eqref{him-p4-e-3.26} and \eqref{him-p4-e-3.26-a}, we obtain
\begin{align} \label{him-p4-e-3.27}
T(1)+K(-1)
=& \int \limits _{0}^{1}\frac{1}{s}\int \limits _{0}^{s}M_{k'}(t) M_{\phi}(t) \, dt\, ds+ K(-1) \\ \nonumber
\geq & \int \limits _{0}^{1}\frac{1}{s}\int \limits _{0}^{s}K'(t) \phi(t) \, dt\, ds + K(-1) \\ \nonumber
=& \int \limits _{0}^{1}\frac{1}{s}\int \limits _{0}^{s}\left(tK''(t)+K'(t)\right) \, dt\, ds +K(-1) \\ \nonumber
=& \int \limits _{0}^{1}\frac{1}{s}\left(sK'(s)-K(s)+K(s)\right)\, ds+ K(-1)\\ \nonumber
=& K(1)+K(-1)>0.
\end{align}
Similarly, by using \eqref{him-p4-e-3.26} and  \eqref{him-p4-e-3.26-a}, we obtain 
\begin{align*}
T'(r)+K'(-r)  &= \frac{1}{r}\int \limits _{0}^{r}M_{K'}(t) M_{\phi}(t) \, \,dt+ K'(-r) \\ 
 &\geq  \frac{1}{r}\int \limits _{0}^{r}K'(t) \phi(t) \, \,dt+ K'(-r) \\ \nonumber
 &= \frac{1}{r}\int \limits _{0}^{r}\left(K'(t)+tK''(t)\right) \, \,dt+ K'(-r) \\ 
  &= K'(r)+K'(-r)>0
\end{align*} 
and hence
\begin{equation} \label{him-p4-e-3.28}
\int\limits_{0}^{1} \left(T'(t)+K'(-t)\right)t \,\, dt>0.
\end{equation}
Combining \eqref{him-p4-e-3.27} and \eqref{him-p4-e-3.28} with \eqref{him-p4-e-3.25}, we obtain $H_{2}(1)>0$.
Similarly, using \eqref{him-p4-e-3.26} and  \eqref{him-p4-e-3.26-a}, we obatin 
$$
H_{2}(0)=-L(1,\alpha)=-K(-1)(1-|\alpha|)+|\alpha|\int \limits_{0}^{1}-K(-t)\,\, dt<0.
$$
Therefore $H_{2}$ has a root in $(0,1)$. Let $r_{f}$ be the smallest root of $H_{2}$ in $(0,1)$. Then $R_{\mathcal{C}_{c}}(r)\leq L(1,\alpha)$ for $r\leq r_{f}$. Combining the inequalities \eqref{him-p4-e-3.13} and \eqref{him-p4-e-3.14} with the fact that $R_{\mathcal{C}_{c}}(r)\leq L(1,\alpha)$ for $r\leq r_{f}$, we obtain 
$$
|z|+\sum_{n=2}^{\infty} (|a_{n}|+|b_{n}|)r^{n} \leq d(f(0), \partial f(\mathbb{D}))
$$
for $|z|=r\leq \min \{1/3, r_{f}\}$. This completes the proof.
\end{pf} 

\noindent\textbf{Acknowledgement:}  The first author thank SERB-MATRICS and the second author thank CSIR for their financial support.

\end{document}